\documentclass[12pt]{article}

\begin{document}

\begin{center}
\textbf{A Generalized Parameter Imbedding Method \\[0pt]
\vspace{12pt} }

\textbf{Cheng Geng \\[0pt]
\vspace{5pt} }

\textbf{Department of Astronomy and Applied Physics, USTC, Hefei, Anhui
230026, P. R. China \\[0pt]
\vspace{12pt} }

\textbf{T. K. Kuo \\[0pt]
\vspace{5pt} }

\textbf{Physics Department, Purdue University, W. Lafayette, IN 47907 \\[0pt]
\vspace{12pt} }

\textbf{ABSTRACT}
\end{center}

\vspace{12pt}

The parameter imbedding method for the Fredholm equation converts it into an
initial value problem in its parameter $\lambda$. We establish the method
for general operator equations of the form $[I + \hat{f}(\lambda)]\psi =
\phi $. It is particularly useful for studying spontaneous symmetry breaking
problems, such as contained in the nonlinear Schwinger-Dyson equation. 
\vspace{12pt}

\pagebreak

\noindent \textbf{I. Introduction} \vspace{8pt}

A very useful technique in the solution of the Fredholm integral equation is
the parameter imbedding method.$^{1}$ Consider the equation 
$$
\psi (x)=\lambda \int ~K(x,y)\psi (y)dy+\phi (x).\eqno{(1)} 
$$

\noindent It is well-known that for $\phi = 0$, the equation has nontrivial
solution only for certain eigenvalues $\lambda = \lambda_i$, and that for $%
\lambda \not{= }\lambda_i$, the solution is given by 
$$
\psi(x) = \phi(x) + \lambda \int~R(x,y,\lambda) \phi(y)dy, \eqno{(2)} 
$$

\noindent where the resolvent function can be written in the form 
$$
R(x,y, \lambda) = D(x,y, \lambda)/d(\lambda). \eqno{(3)} 
$$

\noindent The parameter imbedding method consists of a pair of differential
equations relating $D(x,y,\lambda)$ and $d(\lambda)$,

$$
d^\prime(\lambda) = \frac{d}{d\lambda} d(\lambda) = - \int~D(x,x,\lambda)dx %
\eqno{(4)} 
$$

$$
\frac{d}{d \lambda} D(x,y,\lambda) = D(x,y,\lambda) \frac{d^\prime(\lambda)%
} {d(\lambda)} + \frac{1}{d(\lambda)}~\int D(x,z,\lambda) D(z,y,\lambda)dz, %
\eqno{(5)} 
$$

\noindent with the initial conditions 
$$
d(0) = 1, \eqno{(6)} 
$$

$$
D(x,y,0) = K(x,y). \eqno{(7)} 
$$

\noindent These equations can be integrated numerically with respect to $%
\lambda$, yielding the solution to the original integral equation. Note that
the eigenvalues of the original equation are just the zeroes of $d(\lambda)$%
, also, the eigenfunctions are proportional to $D(x,y^\ast, \lambda)$, for
any fixed value of $y = y^\ast$. The fact that $d(\lambda)$ and $%
D(x,y,\lambda)$ are analytic facilitate the numerical solutions greatly.
Thus, the parameter imbedding method turns the original integral equation,
with variable x and y, into an initial value problem with variable $\lambda$%
. The latter is often much simpler and explicit solutions can be readily
obtained. For a practical application, the reader is referred to the
solutions of the BS equation obtained earlier.$^2$

However, for a number of problems in mathematical physics$^3$, one
encounters integral equations which are not of the Fredholm type. For
instance, a question of considerable current interest concerns systems which
exhibit spontaneous symmetry breaking. They are often described by
nonlinear, Hammerstein equations which are of the form 
$$
\psi(x, \lambda) = \lambda~\int K(x,y) F(y, \psi(y,\lambda))dy, \eqno{(8)} 
$$

\noindent where F is a nonlinear function of $\psi$. The dependence on the
parameter $\lambda$ is thus rather complex. In fact, it typically has
bifurcation solutions which signal the onset of spontaneous symmetry
breaking. One may attempt to solve these equations by linearization and then
apply the parameter imbedding method. However, because of the complicated $%
\lambda$ dependence, the imbedding method for Fredholm equation can not be
used.

The purpose of this paper is to find a generalization of the parameter
imbedding method which can be used in this situation. To this end we have
studied a general class of operator equations in a Hilbert space. They are
of the form 
$$
[I + \hat{f}(\lambda)]\psi = \phi. \eqno{(9)} 
$$

\noindent where we will assume that the operators $\hat{f}(\lambda)$ and $%
\frac{d}{d \lambda} \hat{f}(\lambda)$ are linear, analytic in $\lambda$, and
of trace class$^{4,5}$. (An operator $\hat{A}$ is of trace class if $tr|A| <
\infty)$. Provided we know the inverse operator, $[I + \hat{f}%
(\lambda)]^{-1} $, the solution would be $\psi = [I + \hat{f}%
(\lambda)]^{-1}\phi$. Let us write 
$$
\hat{R}(\lambda) = [I + \hat{f}(\lambda)]^{-1} = \frac{1}{d(\lambda)} \hat{D}%
(\lambda), \eqno{(10)} 
$$

\noindent where 
$$
d(\lambda) = det[I + \hat{f}(\lambda)]. \eqno{(11)} 
$$

\noindent Then, it is found that the functions $d(\lambda)$ and $\hat{D}%
(\lambda)$ are analytic in $\lambda$ and they obey the following pair of
parameter imbedding equations 
$$
\frac{d}{d \lambda} d(\lambda) = Tr[\frac{d}{d \lambda} \hat{f}(\lambda)
\cdot \hat{D}(\lambda)] , \eqno{(12)} 
$$

$$
\frac{d}{d \lambda} \hat{D}(\lambda) = \frac{\hat{D}(\lambda)} {{d} (\lambda)%
}[\frac{d}{d \lambda} d(\lambda) - \frac{d}{d \lambda} \hat{f}(\lambda)
\cdot \hat{D}(\lambda)] . \eqno{(13)} 
$$

\noindent Further, we can determine the initial conditions for $d(\lambda)$
and $\hat{D}(\lambda)$ at an arbitrary point $\lambda= \lambda_o$. Thus, one
can integrate the imbedding equations and arrive at the solution $R(\lambda)$
for any $\lambda$. The imbedding equations, Eqs. (12) and (13), are very
simple, yet rather general. It should be applicable to a wide class of
problems.

We now comment on the background materials to be used in this work. Since we
need to cite results in functional analysis extensively, it is difficult to
make our paper self-contained. For this reason we will follow closely the
four volume treatise ``Methods of Mathematical Physics", by Reed and Simon$%
^4 $, for notations, concepts, and theorems. The only exception is that, in
this paper, all operators will be denoted with a ``hat", as in $\hat{f}%
(\lambda), \hat{D}(\lambda)$, etc.

This paper is organized as follows. The analytic Fredholm theorem in
functional analysis will be presented in Sec. II. In Sec. III and Sec. IV,
we derive the two parameter imbedding equations, which are the central
results of this work. The initial conditions accompanying these equations
are determined in Sec. V, and concluding remarks are offered in Sec. VI. To
make connections between our operator formalism and that of the usual
Fredholm method, we have shown, in Appendix A, a specific realization which
leads to the familiar results. Finally, we present in Appendix B an explicit
formula for the partial trace, which also establishes its relation with the
Plemelj and Smithies formula.

\vspace{8pt}

\noindent \textbf{II. The Fredholm Theorem} \vspace{8pt}

In functional analysis, the Fredholm integral equation has been abstracted
so that it is applicable to the general operator theory in a Hilbert space.
We now quote the analytic Fredholm theorem which will be useful for this
work.$^6$

\noindent Theorem: Let D be an open connected subset of C. Let $\hat{f} = D
\rightarrow \mathcal{L}(\mathcal{H})$ be an analytic operator-value function
such that $\hat{f}(\lambda)$ is compact for each $\lambda \epsilon D$. Then,
either

\noindent A) $[I + \hat{f}(\lambda)]^{-1}$ exists for no $\lambda \epsilon D$%
, or

\noindent B) $[I + \hat{f}(\lambda)]^{-1}$ exists for all $\lambda \epsilon
D \backslash S$, where $S$ is a discrete subset of $D$ (i.e., a set which
has no limit point in $D$). In this case, $[I + \hat{f}(\lambda)]^{-1}$ is
meromorphic in $D$, analytic in $D\backslash S$, the residues at the poles
are finite rank operators, and if $\lambda \epsilon S$ then $-\hat{f}%
(\lambda) \psi = \psi$ has a nonzero solution in $\mathcal{H}$.

>From this theorem, we see that it is only necessary to ascertain that $[I + 
\hat{f}(\lambda)]^{-1}$ exists for some parameter value $\lambda =
\lambda_o\epsilon D$, then the operator $\hat{f}(\lambda)$ would have the
properties described in B). As we will show later, this amounts to the
condition that $d(\lambda_o) \not{= }0$ for some $\lambda = \lambda_o
\epsilon D$. Hereafter, we assume that this condition is always satisfied.

Recall the definition, 
$$
\hat{R}(\lambda) = [I + \hat{f}(\lambda)]^{-1}. \eqno{(14)} 
$$
As long as $\lambda\epsilon D\backslash S$, we have 
$$
\hat{R}(\lambda) + \hat{f}(\lambda) \hat{R}(\lambda) = I, \eqno{(15)} 
$$
$$
\hat{R}(\lambda) + \hat{R}(\lambda) \hat{f}(\lambda) = I . \eqno{(16)} 
$$
Thus, consistency demands that $\hat{f}(\lambda)$ and $\hat{R}(\lambda)$
commute. Note that $\hat{R}(\lambda)$ and $\hat{f}(\lambda)$ are analytic in 
$\lambda\epsilon D\backslash S$ and $\lambda\epsilon D$, respectively. They
and their derivatives are bounded in $D\backslash S$. For later uses we also
note that 
$$
\frac{d}{d\lambda}\hat{R}(\lambda) = - \hat{R}(\lambda) \frac{d}{d\lambda} 
\hat{f}(\lambda) \hat{R}(\lambda), \eqno{(17)} 
$$
which follows directly from Eqs.(15) and (16).

\vspace{8pt}

\noindent \textbf{III. The First Parameter Imbedding Equation} \vspace{8pt}

We now prove our first parameter imbedding equation, Eq.(12). Recall that we
demand both $\hat{f}(\lambda)$ and $\frac{d}{d\lambda}\hat{f}(\lambda)$ to
be of trace class. Since an operator which is of trace class must also be
compact,$^7$ the results of Sec. II are applicable to $\hat{f}(\lambda)$ and 
$\frac{d}{d\lambda} \hat{f}(\lambda)$.

Let's first define the determinant 
$$
d(\lambda) = det (I + \hat{f}(\lambda)) = \sum_{k=o}^\infty Tr [\wedge^k (%
\hat{f}(\lambda))], \eqno{(18)} 
$$
where we used the definition$^4$ 
$$
Tr[\wedge^k\hat{(A)}] = \sum_{1\leq i_1<...<\infty} (e_{i_1} \wedge e_{i_2}
\wedge...\wedge e_{i_k}, \hat{A}e_{i_1} \wedge... \wedge \hat{A} e_{i_k}), %
\eqno{(19)} 
$$
$e_i$ being a basis of the Hilbert space. Also, $Tr[\wedge^o(\hat{A})] = 1$.
The series defining $d(\lambda)$ turns out to be uniformly convergent, which
we now prove.

Let us define the norm of a trace class operator, $\hat{A}\epsilon \mathcal{I%
}_1$, by

$$
||\hat{A}||_{1}=tr|\hat{A}|.\eqno{(20)} 
$$
This gives $\mathcal{I}_{1}$ a norm topology under which it is a Banach
space satisfying the inequality$^{8}$ 
$$
||\wedge ^{k}(\hat{A})||_{1}\leq ||\hat{A}||_{1}^{k}/k!.\eqno{(21)} 
$$
For any $\lambda \epsilon D,\hat{f}(\lambda )$ is of trace class, hence $||%
\hat{f}||_{1}$ is bounded. Let's write $M=sup||\hat{f}(\lambda )||_{1}$, for 
$\lambda \epsilon D$. Then 
$$
||\wedge ^{k}(\hat{f}(\lambda ))||_{1}\leq ||\hat{f}(\lambda
)||_{1}^{k}/k!\leq M^{k}/k!\eqno{(22)} 
$$
It follows that the series for $d(\lambda )$, Eq.(18), is uniformly and
absolutely convergent. We may thus differentiate term by term 
$$
\frac{d}{d\lambda }d(\lambda )=\sum_{k=o}^{\infty }\frac{d}{d\lambda }%
Tr\lbrack \wedge ^{k}(\hat{f}(\lambda ))\rbrack .\eqno{(23)} 
$$
Or, using Eq.(19), and changing orders in the exterior products, we may
write the result in the form 
$$
\everymath={\displaystyle}\frac{d}{d\lambda }d(\lambda )=\sum_{k=o}^{\infty
}k\sum_{1\leq i_{1}<...i_{k}<\infty }(e_{i_{1}}\wedge ...\wedge e_{i_{k}},%
\frac{d}{d\lambda }\hat{f}(\lambda )e_{i_{1}}\wedge ...\wedge \hat{f}%
(\lambda )e_{i_{k}})\eqno{(24)} 
$$
From Appendix B, Eq.(B.1), this may be written as 
$$
\frac{d}{d\lambda }d(\lambda )=\sum_{k=o}^{\infty }kTr\lbrack \frac{d}{%
d\lambda }\hat{f}(\lambda )\cdot Tr_{k-1}(\wedge ^{k}(\hat{f}(\lambda
)))\rbrack ,\eqno{(25)} 
$$
where the partial trace operator is given by Eq.(B.6), 
$$
Tr_{k-1}\lbrack \wedge ^{k}(\hat{f}(\lambda )\rbrack =\sum_{m=1}^{k}\frac{%
(-1)^{m+1}}{k-m+1}\hat{f}^{m-1}(\lambda )Tr~\lbrack \wedge ^{k-m}(\hat{f}%
(\lambda ))\rbrack \eqno{(26)} 
$$
Note that this equation for $\frac{d}{d\lambda }d(\lambda )$ is well-defined
since $\frac{d}{d\lambda }\hat{f}(\lambda )$ is of trace class, $%
Tr_{k-1}\lbrack \wedge ^{k}(\hat{f}(\lambda ))\rbrack $ is bounded or of
trace class, and the product of a bounded operator and one of trace class is
also of trace class.$^{9}$ In addition, the series defined in Eq.(25) is
absolutely and uniformly convergent, following a similar proof for the
series of $d(\lambda )$, Eq.(18). We may thus interchange $Tr$ and $\sum_{k}$%
, resulting in 
$$
\frac{d}{d\lambda }d(\lambda )=Tr\lbrack \frac{d}{d\lambda }\hat{f}(\lambda
)\sum_{k=1}^{\infty }kTr_{k-1}\lbrack \wedge ^{k}(\hat{f}(\lambda ))\rbrack .%
\eqno{(27)} 
$$
Introducing the definition 
$$
\hat{D}(\lambda )=\sum_{k=1}^{\infty }kTr_{k-1}\lbrack \wedge ^{k}(\hat{f}%
(\lambda ))\rbrack ,\eqno{(28)} 
$$
we have 
$$
\frac{d}{d\lambda }d(\lambda )=Tr\lbrack \frac{d}{d\lambda }\hat{f}(\lambda
)\cdot \hat{D}(\lambda )\rbrack .\eqno{(29)} 
$$
This is the first parameter imbedding equation. However, it remains to
establish the relation $\hat{R}(\lambda )=\lbrack I+\hat{f}(\lambda )\rbrack
^{-1}=\hat{D}(\lambda )/d(\lambda )$, which we will do in the next section.

\vspace{8pt}

\noindent \textbf{IV. The Relation $\hat{R}(\lambda) = \hat{D}%
(\lambda)/d(\lambda)$ and the Second Parameter Imbedding Equation}. \vspace{%
8pt}

>From Eqs. (15-16), in order to verify that $\hat{R}(\lambda) = \hat{D}%
(\lambda)/ d(\lambda)$, it is sufficient to establish that 
$$
\hat{D}(\lambda) + \hat{f}(\lambda) \hat{D}(\lambda) = d(\lambda). %
\eqno{(30)} 
$$

$$
\hat{D}(\lambda)+ \hat{D}(\lambda) \hat{f}(\lambda) = d(\lambda). \eqno{(31)}
$$

\noindent Let us first note that $\hat{D}(\lambda)$ and $\hat{f}(\lambda)$
commute. This follows since the partial traces in $\hat{D}(\lambda)$ are all
polynomials in $\hat{f}(\lambda)$. Thus, we need only to examine Eq. (30).
Substituting the formulae for $d(\lambda)$ and $\hat{D}(\lambda)$ from Eqs.
(18) and (28), we have 
$$
\everymath = {\displaystyle} 
\begin{array}{rcl}
& \sum_{k=1}^{\infty}~k~Tr_{k-1} [\wedge^k(\hat{f}(\lambda))] + \hat{f}%
(\lambda)~ \sum_{k=1}^\infty~k~Tr_{k-1} [\wedge^k(\hat{f}(\lambda))] &  \\ 
& = \sum_{k=1}^\infty~Tr[\wedge^{k-1}(\hat{f}(\lambda))]. & 
\end{array}
\eqno{(32)} 
$$

In this equation, we have used the convention $Tr[\wedge^\circ(\hat{f}%
(\lambda))] = 1$, $\hat{f}(\lambda)^\circ = 1$. Eq. (32) is valid since, for
any k, 
$$
\begin{array}{rcl}
& k~Tr_{k-1} [\wedge^k(\hat{f}(\lambda))] + \hat{f}(\lambda)(k-1)
Tr_{k-2}[\wedge^{k-1}(\hat{f}(\lambda))] &  \\ 
& = Tr [\wedge^{k-1}(\hat{f}(\lambda))]. & 
\end{array}
\eqno{(33)} 
$$

\noindent This follows from the Plemelj-Smithies formula, Eq. (B.7), when we
evaluate the determinant by expanding along the first row and using Eq.
(B.9). Note also that the special case, k=1, amounts to $Tr_0(\hat{f}%
(\lambda)) = 1$, which follows from Eq. (B.6).

Having thus verified the relation $\hat{R}(\lambda) = \hat{D}%
(\lambda)/d(\lambda)$, we can now proceed to derive the second parameter
imbedding equation. This we obtain by substituting $\hat{R}(\lambda) = \hat{D%
}(\lambda)/d(\lambda)$ in Eq. (17), derived in Sec. II. It follows
immediately that 
$$
\frac{d}{d\lambda} \hat{D}(\lambda) = \frac{\hat{D}(\lambda)}{d(\lambda)}~ [%
\frac{d}{d \lambda} d(\lambda) - \frac{d}{d \lambda} \hat{f}(\lambda) \hat{D}%
(\lambda)], \eqno{(34)} 
$$

\noindent which is the second imbedding equation. Further, we shall now
deduce that $\hat{D}(\lambda )$ is analytic for all $\lambda \epsilon
D\backslash S$. It is already known that $d(\lambda )$ is analytic for $%
\lambda \epsilon D$. Also, if $d(\lambda )\not{\equiv }0$, the zeroes of $%
d(\lambda )$ are isolated. From the Fredholm theorem in Sec. II, $\hat{R}%
(\lambda )=[I+\hat{f}(\lambda )]^{-1}$ exists for $\lambda \epsilon
D\backslash S$ and is a meromorphic function, whose residues at its poles
are operators of finite rank. Since we already know that $\hat{D}(\lambda )=%
\hat{R}(\lambda )d(\lambda )$ is analytic in $D\backslash S$, to prove its
analyticity in all of D, we need only show that the rank of its zero of $%
d(\lambda )$ is greater than or equal to that of the pole of $\hat{R}%
(\lambda )$. For this purpose we may introduce the operator $\hat{P}%
_{\lambda _{o}}$ which projects $\mathcal{H}$ into the vector space spanned
by the eigenvectors associated with $\lambda _{o}\epsilon S$. This vector
space, $\{\lambda _{o}\}$, is finite dimensional and has dimension $%
dim\{\lambda _{o}\}$. We decompose $\hat{R}(\lambda )$ with respect to $\hat{%
P}_{\lambda _{o}}$ and $(1-\hat{P}_{\lambda _{o})}$, 
$$
\hat{R}(\lambda )=[I+\hat{f}(\lambda )\hat{P}_{\lambda _{o}}]^{-1}\hat{P}%
_{\lambda _{o}}+[I+\hat{f}(\lambda )(1-\hat{P}_{\lambda _{o}})]^{-1}(1-\hat{P%
}_{\lambda _{o})}.\eqno{(35)} 
$$

\noindent At $\lambda =\lambda _{o}$, the second term is nonsingular, while
for the first {term, $\hat{f}(\lambda )\hat{P}_{\lambda _{o}}$ is an
operator in the finite dimensional vector space $\{\lambda _{o}\}$. The
inverse operator $\lbrack I+\hat{f}(\lambda )\hat{P}_{\lambda _{o}}\rbrack
^{-1}$ can thus be obtained with the usual finite-dimensional techniques and
is given by the ratio of a co-factor and the determinant of $\lbrack I+\hat{f%
}(\lambda )\hat{P}_{\lambda _{o}}\rbrack $. The latter is a polynomial of
degree no higher than dim $\{\lambda _{o}\}$, so that the rank of its zero
must be $\leq $ dim $\{\lambda _{o}\}$. Thus, we have proved that $\hat{D}%
(\lambda )=\hat{R}(\lambda )d(\lambda )$ is an analytic function in all of D
(including S). }

\vspace{8pt}

\noindent \textbf{V. The Initial Values of $d(\lambda)$ and $\hat{D}%
(\lambda) $.} \vspace{8pt}

The parameter imbedding equations, given by Eq.(12) and (13), can be
integrated provided we have the initial values of $d(\lambda)$ and $\hat{D}%
(\lambda)$ at some point, say $\lambda = \lambda_o$. This problem becomes
trivial if it happens that $\hat{f}(\lambda_o) = 0$. In this case $\hat{R}%
(\lambda_o) = [I + \hat{f}(\lambda_o)]^{-1} = I$ so that $d(\lambda_o) = 1$
and $\hat{D}(\lambda_o) = I$. In general, we can enforce this condition by
introducing another parameter.

As long as $d(\lambda )\not{\equiv}0$, we can always find a $\lambda _{o}$
so that $d(\lambda _{o})\not{=}0$. Let us consider a new operator equation $%
\lbrack I+\xi \hat{f}(\lambda _{o})\rbrack \psi =\phi $, where $\xi $ is a
complex parameter. The operator-valued functions $\xi \hat{f}(\lambda _{o})$
and $\frac{d}{d\xi }(\xi \hat{f}(\lambda _{o}))=\hat{f}(\lambda _{o})$ are
obviously analytic in $\xi $ and of trace class. Defining $\lbrack I+\xi 
\hat{f}(\lambda _{o})\rbrack ^{-1}=\hat{D}(\xi ,\lambda _{o})/d(\xi ,\lambda
_{o})$, we have the imbedding equations with respect to $\xi $: 
$$
\frac{d}{d\xi }d(\xi ,\lambda _{o})=Tr\lbrack \hat{f}(\lambda _{o})\hat{D}%
(\xi ,\lambda _{o})\rbrack ,\eqno{(36)} 
$$
$$
\frac{d}{d\xi }\hat{D}(\xi ,\lambda _{o})=\frac{\hat{D}(\xi ,\lambda _{o})}{%
d(\xi ,\lambda _{o})}\lbrack \frac{d}{d\xi }d(\xi ,\lambda _{o})-\hat{f}%
(\lambda _{o})\hat{D}(\xi ,\lambda _{o})\rbrack \eqno{(37)} 
$$
For these equations, the initial values at $\xi =0$ are obviously 
$$
d(0,\lambda _{o})=1,\eqno{(38)} 
$$
$$
\hat{D}(0,\lambda _{o})=I.\eqno{(39)} 
$$
We may therefore integrate the imbedding equations to the point $\xi =1$,
where 
$$
d(1,\lambda _{o})=d(\lambda _{o}),\eqno{(40)} 
$$
$$
\hat{D}(1,\lambda _{o})=\hat{D}(\lambda _{o}).\eqno{(41)} 
$$
Provided that $d(\lambda _{o})\not{=}0$, these are then the initial values
we need at $\lambda =\lambda _{o}$, which can be used with the original
Eq.(12) and (13).

\vspace{8pt}

\noindent \textbf{VI. Concluding Remarks} \vspace{8pt}

In this work we have obtained a general method to solve the operator
equation $[I + \hat{f}(\lambda)]\psi = \phi$. If $\hat{f}(\lambda)$ and $%
\frac{d}{d \lambda} \hat{f}(\lambda)$ are linear, analytic in $\lambda$, and
of trace class, then there is a pair of simple differential equations for
the analytic functions $\hat{D}(\lambda)$ and $d(\lambda)$, defined by $\hat{%
R}(\lambda) = [I + \hat{f}(\lambda)]^{-1} = \hat{D}(\lambda)/d(\lambda)$.
The integration of these equations is rather straightforward so that one can
obtain explicitly the functions $d(\lambda)$ and $\hat{D}(\lambda)$.

Since the parameter imbedding method tackles a problem directly with respect
to the parameter $\lambda$, it is particularly useful when the solution
exhibits intriguing properties in $\lambda$. An example is the Hammerstein
equation, which includes the familiar nonlinear Schwinger-Dyson (SD)
equation. There are several advantages in the application of the imbedding
method to this type of problems. First of all,the SD equation is known to
have bifurcation solution. However, the criterion of bifurcation is
contained in the function $d(\lambda)$. As long as $d(\lambda) \not{= }0$,
the implicit function theorem is valid and we have a unique solution. Thus,
at each step of the integration of the parameter imbedding equations, the
obtained value of $d(\lambda)$ enables us to determine whether there is a
unique solution at $\lambda$. Therefore, we can straightforwardly arrive at
a solution for a wide range of the parameter $\lambda$. Another nice feature
of the imbedding method is that, at each step of the integration, the result
obtained is directly the solution in its final form. It is thus very
economical, computationally. Finally, the analyticity properties of $%
d(\lambda)$ and $\hat{D}(\lambda)$ are extremely useful, since we can always
choose contours around possible singularities and integrate over smooth
functions. All of the numerical calculations are therefore routine.

These points are well illustrated in a concrete example, which we present in
a separate paper where we have found numerical solutions to the SD equation.$%
^{10}$ Our imbedding method enables one to obtain explicit bifurcation
solutions to the SD equation. In the process we also uncovered some of its
hitherto unknown solutions.

The imbedding equations discussed in this work are very general. There seems
to be little obstacle in applying them to a number of interesting problems
in mathematical physics. We hope to turn to these questions in the future.

\vspace{8pt}

\noindent \textbf{APPENDIX \ A: }$\mathbf{\hat{R}(\lambda )}$ \textbf{AND
THE FREDHOLM RESOLVENT} \vspace{8pt}

In this appendix, we examine the relation between the operator $\hat{R}%
(\lambda )$ and the usual Fredholm resolvent. They are not identical, which
is why imbedding equations for the two cases also differ. The Fredholm
equation (1) can be written as $\lbrack I+\hat{f}(\lambda )\rbrack \psi
=\phi $ with the definition 
$$
\hat{f}(\lambda )=-\lambda \int ~K(x,y)\psi (y)dy.\eqno{(A1)} 
$$
Its solution \lbrack Eq.(2)\rbrack\ may be written in the form 
$$
\psi =\lbrack I+\lambda R\rbrack \phi \eqno{(A2)} 
$$
with 
$$
R\phi =\int ~R(x,y,\lambda )\phi (y)dy.\eqno{(A3)} 
$$
Thus, the relation between $\hat{R}(\lambda )=\lbrack I+\hat{f}(\lambda
)\rbrack ^{-1}$ and $R$ is 
$$
I+\lambda R=\hat{R}(\lambda )=\lbrack I+\hat{f}(\lambda )\rbrack ^{-1},%
\eqno{(A4)} 
$$
or 
$$
\hat{f}(\lambda )+\lambda R\lbrack I+\hat{f}(\lambda )\rbrack =0.\eqno{(A5)} 
$$
Thus 
$$
R=-\frac{1}{\lambda }\hat{f}(\lambda )\hat{R}(\lambda ).\eqno{(A6)} 
$$

We may now proceed to find a relation of the operator imbedding equations
(12) and (13) for the Fredholm equations. Introducing the coordinate basis $%
\left| x\right\rangle $, so that $\left\langle x\right| \psi \rangle =\psi
(x)$, and defining 
$$
\left\langle x\right| \hat{f}(\lambda )\left| y\right\rangle =-\lambda
K(x,y),\eqno{(A7)} 
$$
$$
\left\langle x\right| \hat{D}(\lambda )\left| y\right\rangle =\hat{D}%
(x,y,\lambda ),\eqno{(A8)} 
$$
then Eq.(12) becomes 
\[
\frac{d}{d\lambda }d(\lambda )=Tr\lbrack \frac{\partial }{\partial \lambda }%
\hat{f}(\lambda )\hat{D}(\lambda )\rbrack =\int \int \left\langle x\right| 
\frac{\partial }{\partial \lambda }\hat{f}(\lambda )\left| y\right\rangle
\left\langle y\right| \hat{D}(\lambda )\left| x\right\rangle dxdy 
\]
$$
=-\int \int K(x,y)\hat{D}(x,y,\lambda )dxdy.\eqno{(A9)} 
$$
Also, taking the matrix element of Eq.(13), we have 
$$
d(\lambda )\frac{\partial }{\partial \lambda }\hat{D}(x,y,\lambda )=\hat{D}%
(x,y,\lambda )\frac{d}{d\lambda }d(\lambda )+\int \hat{D}(x,x^{/},\lambda
)K(x^{/},z)\hat{D}(z,y,\lambda )dx^{/}dz,\eqno{(A10)} 
$$
but Eq.(A6) implies that 
$$
R(x,y,\lambda )=\int K(x,z)\hat{R}(z,y,\lambda )dz,\eqno{(A11)} 
$$
so that 
$$
D(x,y,\lambda )=\int K(x,z)\hat{D}(z,y,\lambda )dz,\eqno{(A12)} 
$$
Thus, Eq.(A9) is precisely Eq.(4). Also, Eq.(A10) is equivalent to Eq.(5) by
multiplying with $K(x^{/},x)$ and integrating with respect to $x$.

Although the two approaches are equivalent, the use of $\hat{R}=\lbrack I+%
\hat{f}(\lambda )\rbrack ^{-1}$ is convenient for general formulation. This
is because the operator $\hat{f}(\lambda )$ is compact, but in an infinite
dimensional space, its inverse is not bounded.$^{11}$ On the other hand, the
identity operator is bounded, but not compact. The combination $\lbrack I+%
\hat{f}(\lambda )\rbrack ^{-1}$ strikes a balance and has nice invertibility
properties similar to those of operators in infinite dimensional spaces. As
for the Fredholm resolvent, we note that from Eq.(A5), it satisfies $%
I+\lambda R\lbrack I+\hat{f}^{-1}\rbrack =0$, or $R=-\lbrack I+\hat{f}%
^{-1}\rbrack ^{-1}/\lambda $. The appearance of the unbounded operator $\hat{%
f}^{-1}$ makes $R$ less useful compared to $\hat{R}$.\textbf{\vspace{1pt}}

\vspace{8pt}

\noindent \textbf{APPENDIX \ B: PARTIAL TRACE} \textbf{AND ITS EXPLICIT FORM}
\vspace{8pt}

Important properties of the partial trace will be examined in this appendix.
For this, we follow Simon,$^{12}$ where the partial trace $Tr_{k-1}$ is
defined as an operator which projects operators in $\mathcal{I}_{1}(\otimes
^{k}\mathcal{H})$to $\mathcal{I}_{1}(\mathcal{H})$. Its satisfies the
relation 
$$
\everymath={\displaystyle}Tr\lbrack \hat{C}Tr_{k-1}(\wedge ^{k}\hat{(A))}%
\rbrack =\sum_{1\leq i_{1}<...i_{k}<\infty }(e_{i_{1}}\wedge ...\wedge
e_{i_{k}},\hat{C}e_{i_{1}}\wedge \hat{A}e_{i_{2}}\wedge ...\wedge \hat{A}%
e_{i_{k}})\eqno{(B1)} 
$$
for any $\hat{C}\in \mathcal{I}_{1}(\mathcal{H})$.

We now proceed to give an explicit realization of the partial trace $%
Tr_{k-1}(\wedge ^{k}\hat{(A))}$. Let us expand the right-hand side of
Eq.(B1) with respect to the first column 
\begin{eqnarray*}
&&\sum_{1\leq i_{1}<...i_{k}<\infty }(e_{i_{1}}\wedge ...\wedge e_{i_{k}},%
\hat{C}e_{i_{1}}\wedge \hat{A}e_{i_{2}}\wedge ...\wedge \hat{A}e_{i_{k}})%
\newline
\\
&=&\frac{1}{k}\sum_{i_{1}=o}^{\infty }k\sum_{1<i_{2}<...<\infty }\lbrack
(e_{i_{1}},\hat{C}e_{i_{1}})(e_{i_{2}}\wedge ...\wedge e_{i_{k}},\hat{A}%
e_{i_{2}}\wedge ...\wedge \hat{A}e_{i_{k}}) \\
&&+\sum_{j=1}^{k-1}(-1)^{j}(e_{i_{j}},\hat{C}e_{i_{1}})(e_{i_{1}}...\wedge
e_{i_{j-1}}\wedge e_{i_{j+1}}\wedge ...,\hat{A}e_{i_{2}}\wedge ...)\rbrack 
\newline
\\
&=&\frac{1}{k}\sum_{i_{1}=o}^{\infty }(e_{i_{1}},\hat{C}e_{i_{1}})Tr\lbrack
\wedge ^{k-1}(\hat{A})\rbrack -\frac{k-1}{k}\sum_{i_{1}=1}^{\infty
}\sum_{1<i_{2}...<\infty }(e_{i_{j}},\hat{C}e_{i_{1}})\newline
\end{eqnarray*}
$$
(...e_{i_{j-1}}\wedge e_{i_{j+1}}...,...\hat{A}e_{i_{j-1}}\wedge \hat{A}%
e_{i_{j+1}}...).\eqno{(B2)}
$$
In deducing the first equality, note that the factor $1/k!$ in the
definition of $\wedge ^{k}(\hat{A})$. For the second equality, note that
there is a factor $(-1)^{j-1}$for bringing $\hat{A}e_{i_{j}}$to the front of
the exterior product. We may now expand the second term in Eq.(B2) into two
terms, etc. Using the relations 
$$
\sum_{i_{1}=1}^{\infty }(e_{i_{1}},\hat{C}e_{i_{1}})=Tr(\hat{C}),\eqno{(B3)}
$$
$$
\sum_{i_{1},i_{j}=1}^{\infty }(e_{i_{j}},\hat{C}e_{i_{1}})(e_{i_{1}},\hat{A}%
e_{i_{j}})=\sum_{i_{j}=1}^{\infty }(e_{i_{j}},\hat{C}\hat{A}e_{i_{j}})=Tr(%
\hat{C}\hat{A}),\eqno{(B4)}
$$
etc., we have 
$$
Tr\lbrack \hat{C}Tr_{k-1}(\wedge ^{k}(\hat{A})\rbrack
=\sum_{m=1}^{k}(-1)^{m+1}\frac{(k-m)!}{k-m+1}Tr(\wedge ^{k-m}(\hat{A}%
))Tr~\lbrack \hat{C}\hat{A}^{m-1}\rbrack .\eqno{(B5)}
$$
It is now clear that the partial trace is given by the explicit formula 
$$
Tr_{k-1}\lbrack \wedge ^{k}(\hat{A})\rbrack =\sum_{m=1}^{k}\frac{(-1)^{m+1}}{%
k-m+1}\hat{A}^{m-1}Tr~\lbrack \wedge ^{k-m}(\hat{A})\rbrack .\eqno{(B6)}
$$
Note that this expression for the partial trace $Tr_{k-1}$ indeed projects $%
\mathcal{L}(\otimes ^{k}\mathcal{H})$ into $\mathcal{L}(\mathcal{H})$,as
required.

Finally, we show that this result for the partial trace is intimately
related to the plemelj-Smithies formula,$^{13}$ which is very important in
the general Fredholm theory. The plemelj-Smithies formula is obtained
through the determinant 
$$
\beta _{k}(\hat{A})=\left| 
\begin{array}{llllll}
\hat{A} & k-1 & 0 & ... & ... & 0 \\ 
\hat{A}^{2} & Tr(\hat{A}) & k-2 & ... & ... & 0 \\ 
... & ... & ... & ... & ... & ... \\ 
\hat{A}^{k} & Tr(\hat{A}^{k-1}) & ... & ... & ... & Tr(\hat{A})
\end{array}
\right| .\eqno{(B7)} 
$$
From the following expression for the trace: 
$$
Tr\lbrack \wedge ^{k}(\hat{A})\rbrack =\frac{1}{k!}\left| 
\begin{array}{llllll}
Tr(\hat{A}) & k-1 & 0 & ... & ... & 0 \\ 
Tr(\hat{A}^{2}) & Tr(\hat{A}) & k-2 & ... & ... & 0 \\ 
... & ... & ... & ... & ... & ... \\ 
Tr(\hat{A}^{k}) & ... & ... & ... & ... & Tr(\hat{A})
\end{array}
\right| ,\eqno{(B8)} 
$$
we can write Eq.(B7) by expanding according to the first column 
$$
\beta _{k}(\hat{A})=\sum_{m=1}^{k}(-1)^{m+1}\frac{k!}{k-m+1}\hat{A}%
^{m}Tr~\lbrack \wedge ^{k-m}(\hat{A})\rbrack =k!\hat{A}Tr_{k-1}\lbrack
\wedge ^{k}(\hat{A})\rbrack .\eqno{(B9)} 
$$
Thus, the explicit formula of the partial trace (B6) is a generalization of
the plemelj-Smithies formula. This result has been crucial in establishing
the parameter imbedding theory for the trace class operator-valued functions.


\begin{center}
\textbf{References}
\end{center}

\vspace{5pt}

\begin{enumerate}
\item  See, e.g., H. Kagiwada and R. Kalaba, Integral Equations via
Imbedding Methods, Addison-Wesley, (1974); also in Solution Method for
Integral Equations, Ed. M.A. Golberg, Plenam Press (1979).

\item  Cheng Geng, Comm. Theor. Phys. \textbf{15}, 219 (1991); \textbf{15},
303 (1991). Cheng Geng and Li Laiyu, ibid, \textbf{15}, 451 (1991).

\item  See, e.g., M. S. Berger, Mathematical Structures of Nonlinear
Science, an Introduction, Kluwer Academic Publishers (1990); Applications of
Nonlinear Analysis in the Physical Science, Eds. H. Amann, N. Bazley and K.
Kirchg\"{a}ssner, Pitman Advanced Publishing Program, (1981).

\item  M. Reed and B. Simon, Methods of Modern Mathematical Physics, Vol. I
- IV, Academic Press (1980).

\item  Ref. 4, Vol. I, p. 207.

\item  Ref. 4, Vol. I, Thm. VI. 14.

\item  Ref. 4, Vol. I, p. 209.

\item  Ref. 4, Vol. IV, p. 323.

\item  Ref. 4, Vol. I, p. 207.

\item  Cheng Geng and T. K. Kuo, Purdue University, Preprint, PURD-TH-94-04.

\item  See, e.g., I. Stakgold, Green's Function and Boundary Value Problem,
John Wiley and Sons, (1979), p. 336.

\item  B. Simon, Adv. Math. \textbf{24}, 244 (1977); Ref. 4, Vol. IV, p. 382.

\item  Ref. 4, Vol. IV, p. 333.
\end{enumerate}

\end{document}